\documentclass[a4paper,11pt,twoside,reqno]{amsart}

\usepackage[utf8]{inputenc}
\usepackage[plainpages=false,pdfpagelabels=true]{hyperref}
\usepackage{amssymb,amsthm}
\usepackage[margin=1in]{geometry}
\usepackage{slashed}

\newtheorem{Satz}{Theorem}[section]
\newtheorem{Prop}[Satz]{Proposition}
\newtheorem{Lem}[Satz]{Lemma}

\theoremstyle{definition}

\newtheorem{Bem}[Satz]{Remark}

\newcommand{\tr}{\operatorname{Tr}}

\newcommand{\dv}{\text{ }dV}
\allowdisplaybreaks[1]

\renewcommand{\epsilon}{\varepsilon}

\newcommand{\R}{\ensuremath{\mathbb{R}}}
\newcommand{\N}{\ensuremath{\mathbb{N}}}

\numberwithin{equation}{section}

\title{On finite energy solutions of 4-harmonic and ES-4-harmonic maps}
\author{Volker Branding}
\date{\today}
\address{University of Vienna, Faculty of Mathematics\\
Oskar-Morgenstern-Platz 1, 1090 Vienna, Austria\\}
\email{volker.branding@univie.ac.at}
\begin{document}

\begin{abstract}
4-harmonic and ES-4-harmonic maps are two generalizations of the well-studied
harmonic map equation which are both given by a nonlinear elliptic partial differential equation
of order eight. Due to the large number of derivatives it is very difficult to 
find any difference in the qualitative behavior of these two variational problems.
In this article we prove that finite energy solutions of both 
4-harmonic and ES-4-harmonic maps from Euclidean space must be trivial.
However, the energy that we require to be finite is different for 4-harmonic and ES-4-harmonic maps
pointing out a first difference between these two variational problems.
\end{abstract} 

\maketitle

\section{Introduction and results}
At the heart of the geometric calculus of variations is the aim to find
interesting maps between Riemannian manifolds. 
This can be achieved by extremizing a given energy functional.

One of the best studied energy functionals for maps between Riemannian manifolds
is the energy of a map \(\phi\colon (M,g)\to (N,h)\) which is 
\begin{align}
\label{energy-map}
E(\phi)=\int_M|d\phi|^2\dv.
\end{align}
The critical points of \eqref{energy-map} are characterized by the vanishing
of the so-called \emph{tension field} which is defined by
\begin{align}
\label{tension-field}
0=\tau(\phi):=\tr_g\bar\nabla d\phi,
\end{align}
where \(\bar\nabla\) represents the connection on \(\phi^\ast TN\).
Solutions of \eqref{tension-field} are called \emph{harmonic maps}
and the latter have been studied intensively in the literature.
The harmonic map equation is a second order semilinear elliptic partial differential
equation. For an overview on the current status of research on harmonic maps
we refer to \cite{MR2389639}.

Recently, many researchers got attracted in energy functionals that contain higher
derivatives extending the energy of a map \eqref{energy-map}.

A possible higher order generalization of harmonic maps is given by the so-called
\emph{polyharmonic maps of order k} or just \emph{k-harmonic maps}. 
These are critical points of the following energy functionals,
where we need to distinguish between polyharmonic maps of even and odd order.
In the even case \((k=2s,s\in\N)\) we set
\begin{align}
\label{poly-energy-even}
E_{2s}(\phi)=\int_M|\bar\Delta^{s-1}\tau(\phi)|^2\dv,
\end{align}
whereas in the odd case \((k=2s+1,s\in\N)\) we have
\begin{align}
\label{poly-energy-odd}
E_{2s+1}(\phi)=\int_M|\bar\nabla\bar\Delta^{s-1}\tau(\phi)|^2\dv.
\end{align}
Here, we use \(\bar\Delta\) to denote the connection Laplacian on the vector bundle \(\phi^\ast TN\).

The first variation of \eqref{poly-energy-even}, \eqref{poly-energy-odd} was calculated in \cite{MR3007953}.
\begin{enumerate}
 \item In the even case (\(k=2s\)) the critical points of \eqref{poly-energy-even} are given by
\begin{align}
\label{tension-2s}
0=\tau_{2s}(\phi):=&\bar\Delta^{2s-1}\tau(\phi)-R^N(\bar\Delta^{2s-2}\tau(\phi),d\phi(e_j))d\phi(e_j) \\
\nonumber&-\sum_{l=1}^{s-1}\bigg(R^N(\bar\nabla_{e_j}\bar\Delta^{s+l-2}\tau(\phi),\bar\Delta^{s-l-1}\tau(\phi))d\phi(e_j) \\
\nonumber&\hspace{1cm}-R^N(\bar\Delta^{s+l-2}\tau(\phi),\bar\nabla_{e_j}\bar\Delta^{s-l-1}\tau(\phi))d\phi(e_j)
\bigg).
\end{align}
\item In the odd case (\(k=2s+1\)) the critical points of \eqref{poly-energy-odd} are given by
\begin{align}
\label{tension-2s+1}
0=\tau_{2s+1}(\phi):=&\bar\Delta^{2s}\tau(\phi)-R^N(\bar\Delta^{2s-1}\tau(\phi),d\phi(e_j))d\phi(e_j)\\
\nonumber&-\sum_{l=1}^{s-1}\bigg(R^N(\bar\nabla_{e_j}\bar\Delta^{s+l-1}\tau(\phi),\bar\Delta^{s-l-1}\tau(\phi))d\phi(e_j) \\
\nonumber&-R^N(\bar\Delta^{s+l-1}\tau(\phi),\bar\nabla_{e_j}\bar\Delta^{s-l-1}\tau(\phi))d\phi(e_j)
\bigg) \\
&\nonumber-R^N(\bar\nabla_{e_j}\bar\Delta^{s-1}\tau(\phi),\bar\Delta^{s-1}\tau(\phi))d\phi(e_j).
\end{align}
\end{enumerate}
Here, we have set \(\bar\Delta^{-1}=0\), \(\{e_j\},j=1,\ldots,m=\dim M\) denotes an orthonormal basis of \(TM\),
and we are applying the Einstein summation convention.

Another possible generalization of harmonic maps,
first suggested by Eells and Sampson in 1964 \cite{MR172310},
can be obtained by studying the critical points of the following energy functional
\begin{align}
\label{k-energy}
E^{ES}_k(\phi):=\int_M|(d+d^\ast)^k\phi|^2\dv=\int_M|(d+d^\ast)^{k-2}\tau(\phi)|^2\dv,\qquad k=1,2,\ldots
\end{align}
For \(k=1\) this energy functional reduces to the energy of a map \eqref{energy-map}.

In the case of \(k=2\), which is also obtained in \eqref{poly-energy-even} for \(s=1\),
we are led to the bienergy \(E_2(\phi)=E^{ES}_2(\phi)\), whose critical points are called \emph{biharmonic maps}.
For an overview on the latter we refer to the recent book \cite{ou2019}.
For biharmonic maps similar classification results as in this article have been obtained in \cite{MR2604617,MR3834926,MR4040175}.

For \(k=3\) we gain the trienergy of a map \(E_3(\phi)=E^{ES}_3(\phi)\), which corresponds to \eqref{poly-energy-odd}
with \(s=1\), and its critical points are called \emph{triharmonic maps}.
For an overview on triharmonic maps we refer to \cite[Section 4]{MR4007262} and references therein.
Triharmonic curves have recently been studied in \cite{montaldo}.

However, for \(k\geq 4\) the energy functional \eqref{k-energy} contains additional curvature terms
and can in general no longer be written in the form \eqref{poly-energy-even}, \eqref{poly-energy-odd}.
An extensive analysis of \eqref{k-energy} and its critical points was carried out recently in \cite{MR4106647}.

This article is devoted to polyharmonic maps of order \(4\) arising either as critical points of \eqref{poly-energy-even}
or \eqref{k-energy}.

The energy functional for 4-harmonic maps (corresponding to \eqref{poly-energy-even} with \(s=2\)) is given by
\begin{align}
\label{energy-4-harmonic}
E_4(\phi)=\int_M|\bar\Delta\tau(\phi)|^2\dv.
\end{align}

The critical points of \eqref{energy-4-harmonic} are characterized by the 
vanishing of the 4-tension field
\begin{align}
\label{4-tension}
0=\tau_{4}(\phi):=&\bar\Delta^{3}\tau(\phi)-R^N(\bar\Delta^{2}\tau(\phi),d\phi(e_j))d\phi(e_j) \\
\nonumber&+R^N(\tau(\phi),\bar\nabla_{e_j}\bar\Delta\tau(\phi))d\phi(e_j) 
-R^N(\bar\nabla_{e_j}\tau(\phi),\bar\Delta\tau(\phi))d\phi(e_j).
\end{align}
Solutions of \eqref{4-tension} are called \emph{4-harmonic maps}.

The energy functional for ES-4-harmonic maps (corresponding to \eqref{k-energy} with \(k=4\)) is given by
\begin{align}
\label{energy-es4-harmonic}
E^{ES}_4(\phi)&=\int_M|(d+d^\ast)^4\phi|^2\dv\\
\nonumber&=\int_M|\bar\Delta\tau(\phi)|^2\dv
+\frac{1}{2}\int_M|R^N(d\phi(e_i),d\phi(e_j))\tau(\phi)|^2\dv.
\end{align}

The first variation of \eqref{energy-es4-harmonic} 
was calculated in \cite[Section 3]{MR4106647} and is characterized by the vanishing 
of the ES-4-tension field \(\tau_4^{ES}(\varphi)\) given by the following expression
\begin{align}
\label{es-4-tension}
\tau_4^{ES}(\phi)=
\tau_4(\phi)+\hat\tau_4(\phi).
\end{align}
Here, \(\tau_4(\phi)\) denotes the 4-tension field \eqref{4-tension} and the quantity \(\hat\tau_4(\phi)\)
is defined by
\begin{align*}
\hat{\tau}_4(\phi)=-\frac{1}{2}\big(2\xi_1+2d^\ast\Omega_1+\bar\Delta\Omega_0+\tr R^N(d\phi(\cdot),\Omega_0)d\phi(\cdot)\big),
\end{align*}
where we have used the following abbreviations
\begin{align}
\label{variables-omega}
\Omega_0&=R^N(d\phi(e_i),d\phi(e_j))R^N(d\phi(e_i),d\phi(e_j))\tau(\phi), \\
\nonumber\Omega_1(X)&=R^N(R^N(d\phi(X),d\phi(e_j))\tau(\phi),\tau(\phi))d\phi(e_j),\\
\nonumber\xi_1&=-(\nabla_{d\phi(e_j)}R^N)(R^N(d\phi(e_i),d\phi(e_j))\tau(\phi),\tau(\phi))d\phi(e_i).
\end{align}

Note that we use a slightly different notation for the \(\xi_1\) term in \eqref{variables-omega} compared to \cite{MR4106647}.

It can be directly seen that both constant and harmonic maps are absolute minimizers of the higher order energy functionals
\eqref{poly-energy-even}, \eqref{poly-energy-odd} and \eqref{k-energy}. In order to understand the 
mathematical structure of these energy functionals 
it seems important to find conditions that force critical points of these functionals to be constant or harmonic maps.

We will prove the following results for finite energy solutions of \eqref{4-tension}
\begin{Satz}
\label{liouville-4-harmonic}
Let \(\phi\colon\R^m\to N\) be a smooth 4-harmonic map, \(m\neq 8\).
Assume that
\begin{align}
\label{thm-finiteness-4-harmonic-ass}
\int_{\R^m}(|d\phi|^2+|\bar\nabla d\phi|^2+|\bar\nabla^2d\phi|^2+|\bar\nabla^3d\phi|^2)\dv<\infty.
\end{align}
If \(m=2\) then \(\phi\) must be harmonic, if \(m>2\) then \(\phi\) must be constant.
\end{Satz}

The second main result of this article is the following result on finite energy solutions of \eqref{es-4-tension}
\begin{Satz}
\label{liouville-ES-4-harmonic}
Let \(\phi\colon\R^m\to N\) be a smooth ES-4-harmonic map, \(m\neq 8\)
and suppose that \(|R^N|_{L^\infty}<\infty\).
Assume that
\begin{align}
\label{thm-finiteness-es4-harmonic-ass}
\int_{\R^m}(|d\phi|^2+|\bar\nabla d\phi|^2+|\bar\nabla^2d\phi|^2+|\bar\nabla^3d\phi|^2+|d\phi|^4|\bar\nabla d\phi|^2+|d\phi|^6)\dv<\infty.
\end{align}
If \(m=2\) then \(\phi\) must be harmonic, if \(m>2\) then \(\phi\) must be constant.
\end{Satz}

\begin{Bem}
\begin{enumerate}
 \item Due to the last two terms in \eqref{thm-finiteness-es4-harmonic-ass} the assumptions of Theorem \ref{liouville-ES-4-harmonic}
are more restrictive than the assumptions of Theorem \ref{liouville-4-harmonic}. This points out a first difference between 
4-harmonic and ES-4-harmonic maps.

\item It should be stressed that the assumptions \eqref{thm-finiteness-4-harmonic-ass}, \eqref{thm-finiteness-es4-harmonic-ass}
are stronger than demanding the finiteness of the 4-energy \eqref{energy-4-harmonic} or the ES-4-energy \eqref{energy-es4-harmonic}.
A similar phenomenon also appears in corresponding results for biharmonic \cite[Theorem 3.4]{MR2604617} and triharmonic \cite[Theorem 4.1]{MR4007262} maps.

\item The terminology ``finite energy solutions'' is usually employed in corresponding results for harmonic maps with finite energy \eqref{energy-map}.
In this paper we use the same wording to denote solutions of \eqref{4-tension} satisfying \eqref{thm-finiteness-4-harmonic-ass} and also
solutions of \eqref{es-4-tension} satisfying \eqref{thm-finiteness-es4-harmonic-ass}.
\end{enumerate}

\end{Bem}

In addition to Theorems \ref{liouville-4-harmonic} and \ref{liouville-ES-4-harmonic} we want to mention 
another result characterizing the behavior of 4-harmonic maps
which is a special case of a structure theorem for polyharmonic maps established in \cite{MR4184658}.
As the proof of this theorem relies on a different method we have to make the additional
assumption that \(N\) has bounded geometry which means that its curvature tensor and
all of its covariant derivatives are bounded.

\begin{Satz}
\label{liouville-4-harmonic-sobolev}
Let \(\phi\colon\R^m\to N\) be a smooth 4-harmonic map, \(m>6\)
and suppose that \(N\) has bounded geometry.
\begin{enumerate}
 \item Suppose that the following condition holds
\begin{align*}
\int_{\R^m}(|d\phi|^m+|\bar\nabla d\phi|^\frac{m}{2}+|\bar\nabla^2d\phi|^\frac{m}{3})\dv<\epsilon
\end{align*}
 for some \(\epsilon>0\) small enough.
 \item In addition, assume that 
 \begin{align*}
 \int_{\R^m}(|\bar\Delta\tau(\phi)|^2+|\bar\nabla\bar\Delta\tau(\phi)|^2+|\bar\nabla^2\bar\Delta\tau(\phi)|^2)\dv<\infty.
 \end{align*}
\end{enumerate}
Then \(\phi\) must be harmonic.
\end{Satz}

This result can also be extended to the case of ES-4-harmonic maps.
At the heart of the proof of Theorem \ref{liouville-4-harmonic-sobolev} is a Sobolev inequality
which is used to control the lower order terms on the right hand side of \eqref{4-tension}.
However, the extension of this technique to ES-4-harmonic maps would require that \(m>10\) and also
\(\tau(\phi)\in W^{4,2}(\R^m,\phi^\ast TN)\) in addition to the assumptions made in Theorem \ref{liouville-ES-4-harmonic}.

We would like to point out that the method of proof used for Theorems \ref{liouville-4-harmonic} and \ref{liouville-ES-4-harmonic} 
only seems to work on Euclidean space as we are making use of a globally defined conformal vector field.
On the other hand, the method of proof used for Theorem \ref{liouville-4-harmonic-sobolev}
is not restricted to \(\R^m\) but works on all manifolds that admit a \emph{Euclidean type Sobolev inequality}.
For more details on the latter see the introduction of \cite{MR4184658} and references therein.

Theorems \ref{liouville-4-harmonic} and \ref{liouville-ES-4-harmonic} make use of the stress-energy tensor.
For harmonic maps this tensor was calculated in \cite{MR655417},
for biharmonic maps it was given in \cite{MR891928} and later systematically derived in \cite{MR2395125}.
For polyharmonic maps the stress-energy tensor was obtained recently in \cite{MR4007262}.

Throughout this article we will use the following notation:
Indices on the domain manifold will be denoted by Latin letters \(i=1,\ldots,m=\dim M\)
and we will employ Greek letters \(\alpha=1,\ldots,n=\dim N\) for indices on the target manifold.
We will use the following sign convention for the \emph{rough Laplacian} acting on sections of $\phi^{\ast}TN$
\begin{align*}
\bar{\Delta}=d^\ast d =-\big(\bar\nabla_{e_i}\bar\nabla_{e_i}-\bar\nabla_{\nabla^M_{e_i}e_i}\big),
\end{align*}
where $\{e_i\},i=1,\ldots m$ is a local orthonormal frame field tangent to $M$.
Moreover, we employ the summation convention and tacitly sum over repeated indices.
We will often write \(\bar\nabla_i\) instead of \(\bar\nabla_{e_i}\).
Throughout this article we make use of the following sign convention for the curvature of a connection
\begin{align*}
R(X,Y)Z=\nabla_X\nabla_YZ-\nabla_Y\nabla_XZ-\nabla_{[X,Y]}Z
\end{align*}
for given vector fields \(X,Y,Z\).
The letter \(C\) will always represent a positive constant whose value may change from line to line.

This article is organized as follows:
In Section 2 we prove Theorem \ref{liouville-4-harmonic}.
Afterwards, in Section 3, we derive the stress-energy tensor for ES-4-harmonic maps
and employ it in Section 4 to prove Theorem \ref{liouville-ES-4-harmonic}.

\section{Proof of Theorem \ref{liouville-4-harmonic}}
In this section we will prove Theorem \ref{liouville-4-harmonic}.
The proof relies on the stress energy-tensor associated with \eqref{energy-4-harmonic}
which can be obtained by varying \eqref{energy-4-harmonic} with respect to the metric on the domain.
This variation was carried out in detail in \cite[Section 2]{MR4007262},
the resulting stress-energy tensor is given by

\begin{align}
\label{4-harmonic-energy-momentum-tensor}
\nonumber S_{4}(X,Y):=&g(X,Y)\bigg(-\frac{1}{2}|\bar\Delta\tau(\phi)|^2-\langle\tau(\phi),\bar\Delta^{2}\tau(\phi)\rangle-\langle d\phi,\bar\nabla\bar\Delta^{2}\tau(\phi)\rangle 
+\langle\bar\nabla\tau(\phi),\bar\nabla\bar\Delta\tau(\phi)\rangle
\bigg) \\
\nonumber&-\langle\bar\nabla_X\tau(\phi),\bar\nabla_Y\bar\Delta\tau(\phi)\rangle-\langle\bar\nabla_Y\bar\Delta\tau(\phi),\bar\nabla_X\bar\Delta\tau(\phi)\rangle\\
&+\langle d\phi(X),\bar\nabla_Y\bar\Delta^{2}\tau(\phi)\rangle
+\langle d\phi(Y),\bar\nabla_X\bar\Delta^{2}\tau(\phi)\rangle.
\end{align}

It was also shown in \cite[Proposition 2.6]{MR4007262} that the stress-energy tensor \eqref{4-harmonic-energy-momentum-tensor}
satisfies the following conservation law:

\begin{Prop}
Let \(\phi\colon M\to N\) be a smooth map.
Then the stress-energy tensor \eqref{4-harmonic-energy-momentum-tensor} satisfies the following conservation law
\begin{align}
\label{conservation-stress-energy-4-harmonic}
\operatorname{div}S_{4}=-\langle\tau_{4}(\phi),d\phi\rangle.
\end{align}
In particular, \(S_4\) is divergence-free whenever \(\phi\) is a 4-harmonic map, that is a solution of \eqref{4-tension}.
\end{Prop}

Now, for \(R>0\) let \(\eta\in C_0^\infty(\mathbb{R})\) be a smooth cut-off function satisfying \(\eta=1\) for \(|z|\leq R\),
\(\eta=0\) for \(|z|\geq 2R\) and \(|\eta^l(z)|\leq\frac{C}{R^l},l=1,\ldots,4\). 
We define the function \(Y(x):=x\eta(r)\in C_0^\infty(\mathbb{R}^m,\mathbb{R}^m)\) with \(r=|x|\).
It follows directly that
\[
\frac{\partial Y_i}{\partial x^j}=\delta_{ij}\eta(r)+\frac{x_i x_j}{r}\eta'(r).
\]
Due to the conservation law \eqref{conservation-stress-energy-4-harmonic} we have
\begin{align*}
0=-\int_{\R^m}\langle Y,\operatorname{div}S_4\rangle\dv=\int_{\R^m}\frac{\partial Y_i}{\partial x^j}S_4(e_i,e_j)\dv.
\end{align*}
By a direct computation we find
\begin{align}
\label{4-harmonic-definition-hr}
\int_{\R^m} S_4(e_i,e_j)\delta_{ij}\eta(r)\dv
=& \int_{\R^m}\eta(r)\big(
-m(\frac{1}{2}|\bar\Delta\tau(\phi)|^2+\langle\tau(\phi),\bar\Delta^2\tau(\phi)\rangle) \\
\nonumber&+(2-m)\langle d\phi,\bar\nabla\bar\Delta^2\tau(\phi)\rangle
+(m-2)\langle\bar\nabla\tau(\phi),\bar\nabla\bar\Delta\tau(\phi)\rangle
\big)\dv \\
\nonumber:=&\sum_{r=1}^4H_r.
\end{align}

As a next step we manipulate the four terms on the right hand side of \eqref{4-harmonic-definition-hr}.
Note that the \(H_1\)-term already has the form that we need.
Hence, we start by manipulating the \(H_2\)-term as follows
\begin{align*}
\int_{\R^m}\eta(r)\langle\tau(\phi),\bar\Delta^2\tau(\phi)\rangle\dv
=&-\int_{\R^m}\big(\eta(r)\big)_{jj}\langle\tau(\phi),\bar\Delta\tau(\phi)\rangle\dv \\
&-2\int_{\R^m}\big(\eta(r)\big)_{j}\langle\bar\nabla_j\tau(\phi),\bar\Delta\tau(\phi)\rangle\dv 
+\int_{\R^m}\eta(r)|\bar\Delta\tau(\phi)|^2\dv.
\end{align*}
Here, and in the following, a subscript \(j\) denotes the derivative with respect to the \(j\)-th coordinate variable
in \(\R^m\).

For the \(H_3\)-term we obtain
\begin{align*}
\int_{\R^m}\eta(r)\langle d\phi,\bar\nabla\bar\Delta^2\tau(\phi)\rangle\dv 
=&-\int_{\R^m}\eta(r)\langle\tau(\phi),\bar\Delta^2\tau(\phi)\rangle\dv 
+\int_{\R^m}\big(\eta(r)\big)_{jkk}\langle d\phi(e_j),\bar\Delta\tau(\phi)\rangle\dv \\
&+2\int_{\R^m}\big(\eta(r)\big)_{jk}\langle\bar\nabla_k d\phi(e_j),\bar\Delta\tau(\phi)\rangle\dv \\
&-\int_{\R^m}\big(\eta(r)\big)_j\langle\bar\Delta d\phi(e_j),\bar\Delta\tau(\phi)\rangle\dv.
\end{align*}

Regarding the \(H_4\)-term we get
\begin{align*}
\int_{\R^m}\eta(r)\langle\bar\nabla\tau(\phi),\bar\nabla\bar\Delta\tau(\phi)\rangle\dv
=\int_{\R^m}\eta(r)|\bar\Delta\tau(\phi)|^2\dv
-\int_{\R^m}\big(\eta(r)\big)_j\langle\bar\nabla_j\tau(\phi),\bar\Delta\tau(\phi)\rangle\dv.
\end{align*}

By another direct computation we obtain
\begin{align}
\label{4-harmonic-definition-jr}
\int_{\R^m} S_4(e_i,e_j)\frac{x_i x_j}{r}\eta'(r)\dv
=& \int_{\R^m}\eta'(r)r\big(
-\frac{1}{2}|\bar\Delta\tau(\phi)|^2-\langle\tau(\phi),\bar\Delta^2\tau(\phi)\rangle) \\
\nonumber&-\langle d\phi,\bar\nabla\bar\Delta^2\tau(\phi)\rangle 
+\langle\bar\nabla\tau(\phi),\bar\nabla\bar\Delta\tau(\phi)\rangle
\big)\dv
\\
\nonumber&+2\int_{\R^m}\eta'(r)\frac{x_i x_j}{r}\langle d\phi(e_i),\bar\nabla_j\bar\Delta^2\tau(\phi)\rangle\dv \\
\nonumber&-2\int_{\R^m}\eta'(r)\frac{x_i x_j}{r}\langle \bar\nabla_i\tau(\phi),\bar\nabla_j\bar\Delta\tau(\phi)\rangle\dv\\
\nonumber:=&\sum_{r=1}^6J_r.
\end{align}

Similar as before,
we will now manipulate the \(J_r\)-terms, \(r=2,\ldots,6\) using integration by parts,
the \(J_1\)-term already has the desired form.
Regarding the \(J_2\) and the \(J_3\)-terms we find
\begin{align*}
J_3=&-J_2-\int_{\R^m}\big(r\eta'(r)\big)_{jkk}\langle d\phi(e_j),\bar\Delta\tau(\phi)\rangle\dv 
-2\int_{\R^m}\big(r\eta'(r)\big)_{jk}\langle\bar\nabla_k d\phi(e_j),\bar\Delta\tau(\phi)\rangle\dv \\
&+\int_{\R^m}\big(r\eta'(r)\big)_{j}\langle\bar\Delta d\phi(e_j),\bar\Delta\tau(\phi)\rangle\dv.
\end{align*}

The \(J_4\)-term can easily manipulated to give
\begin{align*}
J_4=\int_{\R^m}\eta'(r)r|\bar\Delta\tau(\phi)|^2\dv
-\int_{\R^m}\big(r\eta'(r)\big)_j\langle\bar\nabla_j\tau(\phi),\bar\Delta\tau(\phi)\rangle\dv.
\end{align*}

Using integration by parts several times we can express the \(J_5\)-term as 
\begin{align*}
\frac{J_5}{2}
=&\int_{\R^m}\big(\eta'(r)\frac{x_i x_j}{r}\big)_{jkk}\langle d\phi(e_i),\bar\Delta\tau(\phi)\rangle\dv
+2\int_{\R^m}\big(\eta'(r)\frac{x_i x_j}{r}\big)_{jk}\langle\bar\nabla_k d\phi(e_i),\bar\Delta\tau(\phi)\rangle\dv \\
&-\int_{\R^m}\big(\eta'(r)\frac{x_i x_j}{r}\big)_{j}\langle\bar\Delta d\phi(e_i),\bar\Delta\tau(\phi)\rangle\dv
+\int_{\R^m}\big(\eta'(r)\frac{x_i x_j}{r}\big)_{kk}\langle\bar\nabla_j d\phi(e_i),\bar\Delta\tau(\phi)\rangle\dv \\
&+2\int_{\R^m}\big(\eta'(r)\frac{x_i x_j}{r}\big)_{k}\langle\bar\nabla_k\bar\nabla_j d\phi(e_i),\bar\Delta\tau(\phi)\rangle\dv 
-\int_{\R^m}\eta'(r)\frac{x_i x_j}{r}\langle\bar\Delta\bar\nabla_j d\phi(e_i),\bar\Delta\tau(\phi)\rangle\dv.
\end{align*}

Finally, for the \(J_6\)-term we get
\begin{align*}
\frac{J_6}{2}=&\int_{\R^m}\eta'(r)\frac{x_i x_j}{r}\langle \bar\nabla_j\bar\nabla_i\tau(\phi),\bar\Delta\tau(\phi)\rangle\dv
+\int_{\R^m}\big(\eta'(r)\frac{x_ix_j}{r}\big)_j\langle\bar\nabla_i\tau(\phi),\bar\Delta\tau(\phi)\rangle\dv.
\end{align*}

Combining \eqref{4-harmonic-definition-hr} and \eqref{4-harmonic-definition-jr}
and using the identities for \(H_r,r=1,\ldots 4\) and \(J_r,r=1\ldots 6\) we can deduce that

\begin{align}
\label{identity-4-harmonic-a}
(4-&\frac{m}{2})\int_{\R^m}\eta(r)|\bar\Delta\tau(\phi)|^2\dv\\
\nonumber=&2\int_{\R^m}\big(\eta(r)\big)_{jj}\langle\tau(\phi),\bar\Delta\tau(\phi)\rangle\dv 
-(m-6)\int_{\R^m}\big(\eta(r)\big)_{j}\langle\bar\nabla_j\tau(\phi),\bar\Delta\tau(\phi)\rangle\dv \\
\nonumber&+(2-m)\int_{\R^m}\big(\eta(r)\big)_{jkk}\langle d\phi(e_j),\bar\Delta\tau(\phi)\rangle\dv 
-(2-m)\int_{\R^m}\big(\eta(r)\big)_j\langle\bar\Delta d\phi(e_j),\bar\Delta\tau(\phi)\rangle\dv\\
\nonumber&+2(2-m)\int_{\R^m}\big(\eta(r)\big)_{jk}\langle\bar\nabla_k d\phi(e_j),\bar\Delta\tau(\phi)\rangle\dv \\
\nonumber&+\frac{1}{2}\int_{\R^m}\eta'(r)r|\bar\Delta\tau(\phi)|^2\dv
-\int_{\R^m}\big(r\eta'(r)\big)_j\langle\bar\nabla_j\tau(\phi),\bar\Delta\tau(\phi)\rangle\dv \\
\nonumber&-\int_{\R^m}\big(r\eta'(r)\big)_{jkk}\langle d\phi(e_j),\bar\Delta\tau(\phi)\rangle\dv 
-2\int_{\R^m}\big(r\eta'(r)\big)_{jk}\langle\bar\nabla_k d\phi(e_j),\bar\Delta\tau(\phi)\rangle\dv \\
\nonumber&+\int_{\R^m}\big(r\eta'(r)\big)_{j}\langle\bar\Delta d\phi(e_j),\bar\Delta\tau(\phi)\rangle\dv \\
\nonumber&+2\int_{\R^m}\big(\eta'(r)\frac{x_i x_j}{r}\big)_{jkk}\langle d\phi(e_i),\bar\Delta\tau(\phi)\rangle\dv
+4\int_{\R^m}\big(\eta'(r)\frac{x_i x_j}{r}\big)_{jk}\langle\bar\nabla_k d\phi(e_i),\bar\Delta\tau(\phi)\rangle\dv \\
\nonumber&-2\int_{\R^m}\big(\eta'(r)\frac{x_i x_j}{r}\big)_{j}\langle\bar\Delta d\phi(e_i),\bar\Delta\tau(\phi)\rangle\dv
+2\int_{\R^m}\big(\eta'(r)\frac{x_i x_j}{r}\big)_{kk}\langle\bar\nabla_j d\phi(e_i),\bar\Delta\tau(\phi)\rangle\dv \\
\nonumber&+4\int_{\R^m}\big(\eta'(r)\frac{x_i x_j}{r}\big)_{k}\langle\bar\nabla_k\bar\nabla_j d\phi(e_i),\bar\Delta\tau(\phi)\rangle\dv 
-2\int_{\R^m}\eta'(r)\frac{x_i x_j}{r}\langle\bar\Delta\bar\nabla_j d\phi(e_i),\bar\Delta\tau(\phi)\rangle\dv \\
\nonumber&+2\int_{\R^m}\eta'(r)\frac{x_i x_j}{r}\langle \bar\nabla_j\bar\nabla_i\tau(\phi),\bar\Delta\tau(\phi)\rangle\dv
+2\int_{\R^m}\big(\eta'(r)\frac{x_ix_j}{r}\big)_j\langle\bar\nabla_i\tau(\phi),\bar\Delta\tau(\phi)\rangle\dv.
\end{align}

In order to estimate the terms on the right hand side of \eqref{identity-4-harmonic-a} we perform the following 
direct calculations and estimate
\begin{align}
\label{estimate-cutoff-a}
\big(\eta(r)\big)_{jkk}&=\eta'''(r)\frac{x_j}{r}\\
\nonumber&\leq\frac{C}{R^3},\\
\nonumber\big(\eta(r)\big)_{jk}&=\eta''(r)\frac{x_kx_j}{r^2}+\eta'(r)\frac{\delta_{jk}}{r}-\eta'(r)\frac{x_jx_k}{r^3}\\
\nonumber&\leq\frac{C}{R^2},\\
\nonumber\big(\eta(r)\big)_{j}&=\eta'(r)\frac{x_j}{r}\\
\nonumber&\leq\frac{C}{R}.
\end{align}

Similarly, we obtain
\begin{align}
\label{estimate-cutoff-b}
\big(\eta'(r)\frac{x_i x_j}{r}\big)_{jkk}&=\eta^{(4)}(r)x_i+3\eta'''(r)\frac{x_i}{r}
\\
\nonumber&\leq \frac{C}{R^3},\\
\nonumber\big(\eta'(r)\frac{x_i x_j}{r}\big)_{jk}&=\eta'''(r)\frac{x_ix_k}{r}+\eta''(r)\delta_{ik}
+\eta''(r)\frac{x_ix_k}{r^2}+\eta'(r)\frac{\delta_{ik}}{r}-\eta'(r)\frac{x_ix_k}{r^3}\\
\nonumber&\leq\frac{C}{R^2},\\
\nonumber\big(\eta'(r)\frac{x_i x_j}{r}\big)_{kk}&=\eta'''(r)\frac{x_ix_j}{r}
+2\eta''(r)\frac{x_ix_j}{r^2}+2\eta'(r)\frac{\delta_{ij}}{r}-2\eta'(r)\frac{x_ix_j}{r^3}\\
\nonumber&\leq\frac{C}{R^2},\\
\nonumber\big(\eta'(r)\frac{x_i x_j}{r}\big)_{k}&=\eta''(r)\frac{x_ix_jx_k}{r^2}+\eta'(r)\frac{\delta_{ik}x_j}{r}
+\eta'(r)\frac{\delta_{jk}x_i}{r}
-\eta'(r)\frac{x_ix_jx_k}{r^3}\\
\nonumber&\leq\frac{C}{R}.
\end{align}

Inserting \eqref{estimate-cutoff-a} and \eqref{estimate-cutoff-b} into \eqref{identity-4-harmonic-a}
and using Young's inequality multiple times we find
\begin{align*}
\int_{\R^m}\eta(r)|\bar\Delta\tau(\phi)|^2\dv\leq&\frac{C}{|8-m|}\big(\frac{1}{R}+\frac{1}{R^2}+\frac{1}{R^3}\big)
\int_{\R^m}(|d\phi|^2+|\bar\nabla d\phi|^2+|\bar\nabla^2d\phi|^2+|\bar\nabla^3d\phi|^2)\dv \\
&+\frac{C}{|8-m|}\int_{B_{2R}\setminus B_R}|\bar\Delta\tau(\phi)|^2\dv.
\end{align*}
Taking the limit \(R\to\infty\) and using the finiteness assumption \eqref{thm-finiteness-4-harmonic-ass}
the calculation from above yields that \(\bar\Delta\tau(\phi)=0\).

At this point, we employ integration by parts
\begin{align*}
0=-\int_{\R^m}\eta^2\langle\underbrace{\bar\Delta\tau(\phi)}_{=0},\tau(\phi)\rangle\dv=-\int_{\R^m}\eta^2|\bar\nabla\tau(\phi)|^2\dv
-2\int_{\R^m}\eta\nabla\eta\langle\bar\nabla\tau(\phi),\tau(\phi)\rangle\dv
\end{align*}
from which we may deduce that
\begin{align*}
\int_{\R^m}\eta^2|\bar\nabla\tau(\phi)|^2\dv\leq\frac{C}{R^2}\int_{\R^m}|\tau(\phi)|^2\dv\leq\frac{C}{R^2}\int_{\R^m}|\bar\nabla d\phi|^2\dv.
\end{align*}
Again, taking the limit \(R\to\infty\) yields that \(\bar\nabla\tau(\phi)=0\).
Testing \(\bar\nabla\tau(\phi)=0\) with \(-\eta^2d\phi\) and performing the same step as before
we can conclude that \(\tau(\phi)=0\).

Now, the claim that \(\phi\) must be trivial if \(m\neq 2\) follows
from a classical result of Sealey \cite[Corollary 1]{MR654088}
which states that there does not exist a non-constant harmonic map
of finite energy \(E(\phi)=\int_{\R^m}|d\phi|^2\dv\) if the domain
is the Euclidean space \(\R^m,m\geq 3\) with the flat metric.

\section{The stress-energy tensor for ES-4-harmonic maps}
In this section we derive the stress-energy tensor associated with ES-4-harmonic maps
by varying the functional \(E^{ES}_4(\phi)\) with respect to the metric on the domain.
We only have to compute the variation with respect to the metric of the second term in \eqref{energy-es4-harmonic}
as the stress-energy tensor for 4-harmonic maps was already derived in \cite[Section 2]{MR4007262}.
In this section we will allow \((M,g)\) to be an arbitrary Riemannian manifold and do
not restrict to the case \(M=\R^m\).

Throughout this section we set
\begin{align}
\label{variation-metric-domain}
\frac{d}{dt}\big|_{t=0}g_{ij}=\omega_{ij},
\end{align}
where \(\omega_{ij}\) is a smooth symmetric \(2\)-tensor on \(M\).

\begin{Lem}
Let \(\phi\colon M\to N\) be a smooth map and consider a variation of the metric on \(M\) as defined in \eqref{variation-metric-domain}.
Then the following formula holds
\begin{align}
\label{variation-es4-a}
\frac{d}{dt}\big|_{t=0}\frac{1}{2}\int_M&|R^N(d\phi(e_k),d\phi(e_l))\tau(\phi)|^2\dv_{g_t} \\
\nonumber=&\int_M\langle R^N(d\phi(e_k),d\phi(e_l))\frac{d}{dt}\big|_{t=0}\tau(\phi),R^N(d\phi(e_k),d\phi(e_l))\tau(\phi)\rangle\dv_{g} \\
\nonumber &-\int_M\langle R^N(d\phi(e_k),d\phi(e_i))\tau(\phi),R^N(d\phi(e_k),d\phi(e_j))\tau(\phi)\rangle\omega^{ij}\dv_{g}\\
\nonumber &+\frac{1}{4}\int_M|R^N(d\phi(e_k),d\phi(e_l))\tau(\phi)|^2\langle\omega,g\rangle\dv_{g}.
\end{align}
\end{Lem}

\begin{proof}
Recall that the variation of the volume element is given by
\begin{align*}
\frac{d}{dt}|_{t=0}\dv_{g_t}=\frac{1}{2}\langle g,\omega\rangle\dv_g.
\end{align*}
The claim then follows by a direct calculation using that
\begin{align*}
\frac{d}{dt}\big|_{t=0}g^{ij}=-\omega^{ij}.
\end{align*}

\end{proof}

To proceed we recall the following lemma (see for example \cite[Lemma 2.2]{MR4007262})
\begin{Lem}
\label{lem-variation-tension}
Let \(\phi\colon M\to N\) be a smooth map and consider a variation of the metric on \(M\) as defined in \eqref{variation-metric-domain}.
The variation of the tension field with respect to the metric on the domain is given by
\begin{align}
\label{variation-tension-field}
\frac{d}{dt}\big|_{t=0}\tau^\alpha(\phi)=-\omega^{ij}(\bar\nabla d\phi)^\alpha_{ij}-(\nabla_i\omega^{ki})d\phi^\alpha(e_k)
+\frac{1}{2}(\nabla^k\tr\omega)d\phi^\alpha(e_k),
\end{align}
where \(\alpha=1,\ldots,n\).
\end{Lem}

This allows us to perform the following computation:

\begin{Lem}
Let \(\phi\colon M\to N\) be a smooth map and consider a variation of the metric on \(M\) as defined in \eqref{variation-metric-domain}.
Then the following formula holds
\begin{align}
\label{variation-es4-b}
\int_M\langle& R^N(d\phi(e_k),d\phi(e_l))\frac{d}{dt}\big|_{t=0}\tau(\phi),R^N(d\phi(e_k),d\phi(e_l))\tau(\phi)\rangle\dv_g \\
\nonumber=&-\int_M\langle \bar\nabla_i\big(R^N(d\phi(e_j),d\phi(e_l))R^N(d\phi(e_j),d\phi(e_l))\tau(\phi)\big),d\phi(e_k)\rangle\omega^{ki}\dv_g \\
\nonumber&+\frac{1}{2}\int_M\langle\bar\nabla^k\big( R^N(d\phi(e_j),d\phi(e_l))R^N(d\phi(e_j),d\phi(e_l))\tau(\phi)\big),d\phi(e_k)\rangle\langle\omega,g\rangle\dv_g \\
\nonumber&-\frac{1}{2}\int_M|R^N(d\phi(e_k),d\phi(e_l))\tau(\phi)|^2\langle\omega,g\rangle\dv_g.
\end{align}

\end{Lem}

\begin{proof}
Using \eqref{variation-tension-field} in the first term on the right hand side of \eqref{variation-es4-a} we obtain
\begin{align*}
\int_M\langle R^N(d\phi(e_k),&d\phi(e_l))\frac{d}{dt}\big|_{t=0}\tau(\phi),R^N(d\phi(e_k),d\phi(e_l))\tau(\phi)\rangle\dv_g \\
=&-\int_M\langle R^N(d\phi(e_k),d\phi(e_l))\bar\nabla_{i}d\phi(e_j),R^N(d\phi(e_k),d\phi(e_l))\tau(\phi)\rangle\omega^{ij}\dv_g\\
&-\int_M\langle R^N(d\phi(e_j),d\phi(e_l))d\phi(e_k),R^N(d\phi(e_j),d\phi(e_l))\tau(\phi)\rangle\nabla_i\omega^{ki}\dv_g\\
&+\frac{1}{2}\int_M\langle R^N(d\phi(e_j),d\phi(e_l))d\phi(e_k),R^N(d\phi(e_j),d\phi(e_l))\tau(\phi)\rangle\nabla^k\tr\omega\dv_g.
\end{align*}
The first two terms can be manipulated as follows
\begin{align*}
-\int_M\langle& R^N(d\phi(e_k),d\phi(e_l))\bar\nabla_{i}d\phi(e_j),R^N(d\phi(e_k),d\phi(e_l))\tau(\phi)\rangle\omega^{ij}\dv_g\\
-\int_M\langle& R^N(d\phi(e_j),d\phi(e_l))d\phi(e_k),R^N(d\phi(e_j),d\phi(e_l))\tau(\phi)\rangle\nabla_i\omega^{ki}\dv_g\\
=&\int_M\langle R^N(d\phi(e_k),d\phi(e_l))R^N(d\phi(e_k),d\phi(e_l))\tau(\phi),\bar\nabla_{i}d\phi(e_j)\rangle\omega^{ij}\dv_g\\
&+\int_M\langle R^N(d\phi(e_j),d\phi(e_l))R^N(d\phi(e_j),d\phi(e_l))\tau(\phi),d\phi(e_k)\rangle\nabla_i\omega^{ki}\dv_g\\
=&-\int_M\langle \bar\nabla_i\big(R^N(d\phi(e_j),d\phi(e_l))R^N(d\phi(e_j),d\phi(e_l))\tau(\phi)\big),d\phi(e_k)\rangle\omega^{ki}\dv_g,
\end{align*}
where we first used the symmetries of the Riemann curvature tensor and applied integration by parts in the second step.

Regarding the third term a similar manipulation yields
\begin{align*}
\int_M\langle& R^N(d\phi(e_j),d\phi(e_l))d\phi(e_k),R^N(d\phi(e_j),d\phi(e_l))\tau(\phi)\rangle\nabla^k\tr\omega\dv_g \\
=&\int_M\langle\bar\nabla^k\big( R^N(d\phi(e_j),d\phi(e_l))R^N(d\phi(e_j),d\phi(e_l))\tau(\phi)\big),d\phi(e_k)\rangle\langle\omega,g\rangle\dv_g \\
&-\int_M|R^N(d\phi(e_k),d\phi(e_l))\tau(\phi)|^2\langle\omega,g\rangle\dv_g.
\end{align*}
The claim then follows from combining the equations.

\end{proof}

We may now give the following
\begin{Prop}
Let \(\phi\colon M\to N\) be a smooth map and consider a variation of the metric on \(M\) as defined in \eqref{variation-metric-domain}.
Then the following formula holds
\begin{align*}
\frac{d}{dt}\big|_{t=0}\frac{1}{2}\int_M&|R^N(d\phi(e_k),d\phi(e_l))\tau(\phi)|^2\dv_{g_t}=\int_M\langle \hat S_4,\omega\rangle\dv_g,
\end{align*}
where the symmetric tensor \(\hat S_4\) is given by
\begin{align}
\label{energy-momentum-e4hat}
\hat S_4(X,Y): =&
-\langle R^N(d\phi(e_k),d\phi(X))\tau(\phi),R^N(d\phi(e_k),d\phi(Y))\tau(\phi)\rangle \\
\nonumber&-\frac{1}{4}|R^N(d\phi(e_k),d\phi(e_l))\tau(\phi)|^2g(X,Y)\\
\nonumber&-\frac{1}{2}\langle\bar\nabla_X\big( R^N(d\phi(e_k),d\phi(e_l))R^N(d\phi(e_k),d\phi(e_l))\tau(\phi)\big),d\phi(Y)\rangle\\
\nonumber&-\frac{1}{2}\langle\bar\nabla_Y\big( R^N(d\phi(e_k),d\phi(e_l))R^N(d\phi(e_k),d\phi(e_l))\tau(\phi)\big),d\phi(X)\rangle\\
\nonumber&+\frac{1}{2}\langle\bar\nabla^k\big( R^N(d\phi(e_k),d\phi(e_l))R^N(d\phi(e_k),d\phi(e_l))\tau(\phi)\big),d\phi(e_k)\rangle g(X,Y).
\end{align}
Here, \(X,Y\) are vector fields on \(M\).
\end{Prop}
\begin{proof}
This follows by combining \eqref{variation-es4-a}, \eqref{variation-es4-b}
and symmetrizing the first term on the right hand side of \eqref{variation-es4-b}.
\end{proof}

\begin{Bem}
In terms of the variables \(\Omega_0,\Omega_1\) defined in \eqref{variables-omega} we may express \eqref{energy-momentum-e4hat} as follows
\begin{align}
\label{energy-momentum-e4hat-omega-variables}
\hat S_4(X,Y)=&
-\langle\Omega_1(Y),d\phi(X)\rangle
+\frac{1}{4}\langle\Omega_0,\tau(\phi)\rangle g(X,Y) \\
\nonumber &+\frac{1}{2}\langle\bar\nabla^k\Omega_0,d\phi(e_k)\rangle g(X,Y)
-\frac{1}{2}\langle\bar\nabla_X\Omega_0,d\phi(Y)\rangle-\frac{1}{2}\langle\bar\nabla_Y\Omega_0,d\phi(X)\rangle.
\end{align}
\end{Bem}

The trace of \eqref{energy-momentum-e4hat} can easily be computed and yields
\begin{align*}
\tr\hat S_4=&(-1-\frac{m}{4})|R^N(d\phi(e_i),d\phi(e_j))\tau(\phi)|^2 \\
\nonumber &+(-1+\frac{m}{2})\langle\bar\nabla^k\big( R^N(d\phi(e_j),d\phi(e_l))R^N(d\phi(e_j),d\phi(e_l))\tau(\phi)\big),d\phi(e_k)\rangle.
\end{align*}

\begin{Bem}
In the case of \(M\) being compact we may use integration by parts to deduce
\begin{align*}
\int_M\tr\hat S_4\dv=\int_M(\frac{m}{4}-2)|R^N(d\phi(e_i),d\phi(e_j))\tau(\phi)|^2\dv.
\end{align*}
This reflects the fact that ES-4-harmonic maps are critical if \(\dim M=8\).
\end{Bem}

Having calculated the variation of \(E^{ES}_4(\phi)\) with respect to the metric
on the domain we may now define the stress-energy tensor for
ES-4-harmonic maps as follows:

\begin{align}
\label{energy-momentum-e4es}
S^{ES}_4(X,Y):=S_4(X,Y)+\hat S_4(X,Y).
\end{align}

The stress-energy tensor \eqref{energy-momentum-e4es} satisfies the following conservation law:

\begin{Satz}
Let \(\phi\colon M\to N\) be a smooth ES-4-harmonic map, that is a smooth solution of 
\begin{align*}
\tau_4^{ES}(\phi)=0,
\end{align*}
where \(\tau_4^{ES}(\phi)\) is defined in \eqref{es-4-tension}.
Then the stress-energy tensor defined in \eqref{energy-momentum-e4es} satisfies
\begin{align}
\label{conservation-stress-energy-ES-4-harmonic}
\operatorname{div}S^{ES}_{4}=-\langle\tau^{ES}_{4}(\phi),d\phi\rangle.
\end{align}
In particular, the stress-energy tensor is divergence-free whenever \(\phi\) is a solution of \(\tau_4^{ES}(\phi)=0\).
\end{Satz}

\begin{proof}
We choose a local orthonormal basis \(e_i,i,\ldots,m\) around a point \(p\in M\) that satisfies \(\nabla_ke_r=0\), where \(1\leq r,k\leq m\).
Then, we calculate
\begin{align}
\label{conservation-stress-energy-identity-aa}
\nabla^j S^4(e_i,e_j)=&-\langle \tau_4(\phi),d\phi(e_i)\rangle \\
\nonumber=&-\langle\xi_1,d\phi(e_i)\rangle
-\langle d^\ast\Omega_1,d\phi(e_i)\rangle
-\frac{1}{2}\langle\bar\Delta\Omega_0,d\phi(e_i)\rangle \\
\nonumber&-\frac{1}{2}\langle\tr R^N(d\phi(\cdot),\Omega_0)d\phi(\cdot),d\phi(e_i)\rangle.
\end{align}
In the first step we made use of \eqref{conservation-stress-energy-4-harmonic} and in the second step
we used that \(\phi\) solves \(\tau_4^{ES}(\phi)=0\).

Now, we will show that the right hand side of \eqref{conservation-stress-energy-identity-aa}
is equal to the negative divergence of \(\hat S^4\).
To this end we calculate using \eqref{energy-momentum-e4hat-omega-variables}
\begin{align}
\label{conservation-stress-energy-identity-ab}
\nonumber\nabla^j\hat S_4(e_i,e_j)=&-\langle\bar\nabla^j\Omega_1(e_j),d\phi(e_i)\rangle
-\langle\Omega_1(e_j),\bar\nabla_jd\phi(e_i)\rangle
-\frac{1}{4}\langle\bar\nabla_i\Omega_0,\tau(\phi)\rangle
+\frac{1}{4}\langle\Omega_0,\bar\nabla_i\tau(\phi)\rangle\\
&+\frac{1}{2}\langle\bar\Delta\Omega_0,d\phi(e_i)\rangle
+\frac{1}{2}\langle R^N(d\phi(e_i),d\phi(e_k))\Omega_0,d\phi(e_k)\rangle.
\end{align}

Using the definition of \(\Omega_1\) given in \eqref{variables-omega} it is easy to derive
\begin{align}
\label{conservation-stress-energy-identity-a}
\langle\Omega_1(e_j),\bar\nabla_jd\phi(e_i)\rangle=
-\langle R^N(d\phi(e_k),\bar\nabla_jd\phi(e_i))\tau(\phi),R^N(d\phi(e_j),d\phi(e_k))\tau(\phi)\rangle.
\end{align}

Moreover, a direct calculation yields
\begin{align}
\label{conservation-stress-energy-identity-b}
\langle\Omega_0,\bar\nabla_i\tau(\phi)\rangle
=-\langle R^N(d\phi(e_k),d\phi(e_j))\bar\nabla_i\tau(\phi),R^N(d\phi(e_k),d\phi(e_j))\tau(\phi)\rangle.
\end{align}

In addition, we find by a direct calculation
\begin{align}
\label{conservation-stress-energy-identity-c}
\langle\bar\nabla_i\Omega_0,\tau(\phi)\rangle=&
2\langle(\nabla_{d\phi(e_i)}R^N)(d\phi(e_k),d\phi(e_j))R^N(d\phi(e_k),d\phi(e_j))\tau(\phi),\tau(\phi)\rangle \\
\nonumber&-4\langle R^N(\bar\nabla_i d\phi(e_k),d\phi(e_l))\tau(\phi),R^N(d\phi(e_k),d\phi(e_l))\tau(\phi)\rangle \\
\nonumber&-\langle R^N(d\phi(e_k),d\phi(e_l))\bar\nabla_i \tau(\phi),R^N(d\phi(e_k),d\phi(e_l))\tau(\phi)\rangle.
\end{align}

Moreover, we manipulate
\begin{align*}
\langle(\nabla&_{d\phi(e_i)}R^N)(d\phi(e_k),d\phi(e_j))R^N(d\phi(e_k),d\phi(e_j))\tau(\phi),\tau(\phi)\rangle \\
=&-\langle(\nabla_{d\phi(e_j)}R^N)(d\phi(e_i),d\phi(e_k))R^N(d\phi(e_k),d\phi(e_j))\tau(\phi),\tau(\phi)\rangle \\
&-\langle(\nabla_{d\phi(e_k)}R^N)(d\phi(e_j),d\phi(e_i))R^N(d\phi(e_k),d\phi(e_j))\tau(\phi),\tau(\phi)\rangle \\
=&-2\langle(\nabla_{d\phi(e_j)}R^N)(d\phi(e_i),d\phi(e_k))R^N(d\phi(e_k),d\phi(e_j))\tau(\phi),\tau(\phi)\rangle \\
=&-2\langle(\nabla_{d\phi(e_j)}R^N)(\tau(\phi),R^N(d\phi(e_k),d\phi(e_j))\tau(\phi))d\phi(e_k),d\phi(e_i)\rangle,
\end{align*}
where we first used the second Bianchi identity and afterwards the symmetries of the Riemannian curvature tensor
in the second and third step.

Combining \eqref{conservation-stress-energy-identity-a}, \eqref{conservation-stress-energy-identity-b} and
\eqref{conservation-stress-energy-identity-c} we get
\begin{align*}
-\langle\Omega_1(e_j),\bar\nabla_jd\phi(e_i)\rangle
-\frac{1}{4}\langle\bar\nabla_i\Omega_0,\tau(\phi)\rangle
+\frac{1}{4}\langle\Omega_0,\bar\nabla_i\tau(\phi)\rangle
=\langle\xi_1,d\phi(e_i)\rangle,
\end{align*}
and together with \eqref{conservation-stress-energy-identity-ab} this completes the proof.
\end{proof}

\begin{Bem}
It was to be expected that the stress-energy tensor associated with the ES-4-energy \eqref{energy-es4-harmonic} is divergence
free. The energy functional \eqref{energy-es4-harmonic} is invariant 
under diffeomorphisms on the domain \(u\colon M\to M\) in the following sense
\begin{align*}
E_4^{ES}(\phi\circ u,u^\ast g)=E_4^{ES}(\phi,g).
\end{align*}
This can be explicitly checked with the methods presented in \cite[Section 2.3]{MR4007262}.
Via Noether's theorem the invariance of the energy functional \eqref{energy-es4-harmonic}
leads to a conserved quantity which is precisely the stress-energy tensor \eqref{energy-momentum-e4es}.
\end{Bem}

\section{Proof of Theorem \ref{liouville-ES-4-harmonic}}
In this section we will prove Theorem \ref{liouville-ES-4-harmonic}.
Our method of proof will be the same as in the proof of Theorem \ref{liouville-4-harmonic}
but instead of the stress-energy tensor for 4-harmonic maps \eqref{energy-4-harmonic} we will now make
use of the stress-energy tensor for ES-4-harmonic maps given by \eqref{energy-momentum-e4es}.

As the stress-energy tensor for ES-4-harmonic maps consists of the stress-energy tensor
for 4-harmonic maps and an additional piece arising from the curvature term \(\hat E^{ES}_4(\phi)\),
which is given by \(\hat S^4\),
we will only have to deal with \(\hat S^4\) as the rest of the calculation will be identical
to the one for 4-harmonic maps.

Choosing the same cutoff function as in the proof of Theorem \ref{liouville-4-harmonic} and
due to the conservation law \eqref{conservation-stress-energy-ES-4-harmonic} we have
\begin{align*}
0=-\int_{\R^m}\langle x\eta(r),\operatorname{div}S^{ES}_4\rangle\dv=\int_{\R^m}\frac{\partial }{\partial x^j}\big(x^i\eta(r)\big)S^{ES}_4(e_i,e_j)\dv.
\end{align*}

Inserting the second term from \eqref{energy-momentum-e4es} into the above equation we find
\begin{align*}
\int_{\R^m}& \hat S_4(e_i,e_j)\delta_{ij}\eta(r)\dv \\
=& \int_{\R^m}\eta(r)\big(
(-1-\frac{m}{4})|R^N(d\phi(e_k),d\phi(e_l))\tau(\phi)|^2\dv \\
&+(-1+\frac{m}{2})\langle\bar\nabla^k\big( R^N(d\phi(e_j),d\phi(e_l))R^N(d\phi(e_j),d\phi(e_l))\tau(\phi)\big),d\phi(e_k)\rangle
\big)\dv.
\end{align*}

Using integration by parts we find
\begin{align*}
\int_{\R^m}&\eta(r)\langle\bar\nabla^k\big( R^N(d\phi(e_j),d\phi(e_l))R^N(d\phi(e_j),d\phi(e_l))\tau(\phi)\big),d\phi(e_k)\rangle\dv \\
=&-\int_{\R^m}\big(\eta(r)\big)_k\langle R^N(d\phi(e_j),d\phi(e_l))R^N(d\phi(e_j),d\phi(e_l))\tau(\phi),d\phi(e_k)\rangle\dv \\
&+\int_{\R^m}\eta(r)|R^N(d\phi(e_k),d\phi(e_l))\tau(\phi)|^2\dv.
\end{align*}

Consequently, we obtain
\begin{align}
\label{identity-ES-4-harmonic-a}
\int_{\R^m}&\hat S_4(e_i,e_j)\delta_{ij}\eta(r)\dv \\
\nonumber=&(-2+\frac{m}{4})\int_{\R^m}\eta(r)|R^N(d\phi(e_k),d\phi(e_l))\tau(\phi)|^2 \dv\\
\nonumber&+(1-\frac{m}{2})\int_{\R^m}\big(\eta(r)\big)_k\langle R^N(d\phi(e_j),d\phi(e_l))R^N(d\phi(e_j),d\phi(e_l))\tau(\phi),d\phi(e_k)\rangle\dv.
\end{align}

Moreover, we find
\begin{align*}
\int_{\R^m}& \hat S_4(e_i,e_j)\frac{x_i x_j}{r}\eta'(r)\dv \\
=&-\frac{1}{4}\int_{\R^m}\eta'(r)r|R^N(d\phi(e_i),d\phi(e_j))\tau(\phi)|^2\dv\\
&+\frac{1}{2}\int_{\R^m}\eta'(r)r\langle\bar\nabla^k\big( R^N(d\phi(e_j),d\phi(e_l))R^N(d\phi(e_j),d\phi(e_l))\tau(\phi)\big),d\phi(e_k)\rangle\dv\\
&-\int_{\R^m}\eta'(r)\frac{x_i x_j}{r}\langle R^N(d\phi(e_k),d\phi(e_i))\tau(\phi),R^N(d\phi(e_k),d\phi(e_j))\tau(\phi)\rangle\dv \\
&-\int_{\R^m}\eta'(r)\frac{x_i x_j}{r}\langle\bar\nabla_i\big( R^N(d\phi(e_k),d\phi(e_l))R^N(d\phi(e_k),d\phi(e_l))\tau(\phi)\big),d\phi(e_j)\rangle\dv.
\end{align*}

In addition, it is straightforward to manipulate
\begin{align*}
\int_{\R^m}&\eta'(r)r\langle\bar\nabla^k\big( R^N(d\phi(e_j),d\phi(e_l))R^N(d\phi(e_j),d\phi(e_l))\tau(\phi)\big),d\phi(e_k)\rangle\dv\\
=&-\int_{\R^m}\big(\eta'(r)r\big)_k\langle R^N(d\phi(e_j),d\phi(e_l))R^N(d\phi(e_j),d\phi(e_l))\tau(\phi),d\phi(e_k)\rangle\dv \\
&+\int_{\R^m}\eta'(r)r|R^N(d\phi(e_i),d\phi(e_j))\tau(\phi)|^2 \dv
\end{align*}
and also
\begin{align*}
\int_{\R^m}&\eta'(r)\frac{x_i x_j}{r}\langle\bar\nabla_i\big( R^N(d\phi(e_k),d\phi(e_l))R^N(d\phi(e_k),d\phi(e_l))\tau(\phi)\big),d\phi(e_j)\rangle\dv \\
=&-\int_{\R^m}\big(\eta'(r)\frac{x_i x_j}{r}\big)_i\langle R^N(d\phi(e_k),d\phi(e_l))R^N(d\phi(e_k),d\phi(e_l))\tau(\phi),d\phi(e_j)\rangle\dv\\
&+\int_{\R^m}\eta'(r)\frac{x_i x_j}{r}\langle R^N(d\phi(e_k),d\phi(e_l))\bar\nabla_id\phi(e_j),R^N(d\phi(e_k),d\phi(e_l))\tau(\phi)\rangle\dv.
\end{align*}

Consequently, we get
\begin{align}
\label{identity-ES-4-harmonic-b}
\int_{\R^m}\hat S_4(e_i,e_j)&\frac{x_i x_j}{r}\eta'(r)\dv=
\frac{1}{4}\int_{\R^m}\eta'(r)r|R^N(d\phi(e_i),d\phi(e_j))\tau(\phi)|^2\dv\\
\nonumber&-\frac{1}{2}\int_{\R^m}\big(\eta'(r)r\big)_k\langle R^N(d\phi(e_j),d\phi(e_l))R^N(d\phi(e_j),d\phi(e_l))\tau(\phi),d\phi(e_k)\rangle\dv \\
\nonumber&-\int_{\R^m}\eta'(r)\frac{x_i x_j}{r}\langle R^N(d\phi(e_k),d\phi(e_i))\tau(\phi),R^N(d\phi(e_k),d\phi(e_j))\tau(\phi)\rangle\dv \\
\nonumber&+\int_{\R^m}\big(\eta'(r)\frac{x_i x_j}{r}\big)_i\langle R^N(d\phi(e_k),d\phi(e_l))R^N(d\phi(e_k),d\phi(e_l))\tau(\phi),d\phi(e_j)\rangle\dv\\
\nonumber&-\int_{\R^m}\eta'(r)\frac{x_i x_j}{r}\langle R^N(d\phi(e_k),d\phi(e_l))\bar\nabla_id\phi(e_k),R^N(d\phi(e_j),d\phi(e_l))\tau(\phi)\rangle\dv.
\end{align}

In order to estimate the terms in \eqref{identity-ES-4-harmonic-a} and \eqref{identity-ES-4-harmonic-b} we use
Young's inequality in the following form
\begin{align*}
\langle R^N(d\phi(e_j),d\phi(e_l))R^N(d\phi(e_j),d\phi(e_l))\tau(\phi),d\phi(e_k)\rangle&\leq
C|d\phi|^5|\bar\nabla d\phi| \\
&\leq C(|d\phi|^6+|d\phi|^4|\bar\nabla d\phi|^2).
\end{align*}

Together with the estimates on the cutoff function \eqref{estimate-cutoff-a}, \eqref{estimate-cutoff-b}
and the estimates obtained in the proof of Theorem \ref{liouville-4-harmonic},
which are \eqref{4-harmonic-definition-hr}, \eqref{4-harmonic-definition-jr} and \eqref{identity-4-harmonic-a},
we get the following inequality
\begin{align*}
\int_{\R^m}&\eta(r)(|\bar\Delta\tau(\phi)|^2\dv+\frac{1}{2}|R^N(d\phi(e_k),d\phi(e_l))\tau(\phi)|^2)\dv \\
\leq &\frac{C}{|8-m|}\big(\frac{1}{R}+\frac{1}{R^2}+\frac{1}{R^3}\big)
\int_{\R^m}(|d\phi|^2+|\bar\nabla d\phi|^2+|\bar\nabla^2d\phi|^2+|\bar\nabla^3d\phi|^2)\dv \\
&+\frac{C}{|8-m|}\frac{1}{R}
\int_{\R^m}(|\bar\nabla d\phi|^2|d\phi|^4+|d\phi|^6)\dv 
+\frac{C}{|8-m|}\int_{B_{2R}\setminus B_R}(|\bar\Delta\tau(\phi)|^2+|d\phi|^4|\tau(\phi)|^2)\dv.
\end{align*}
As long as \(\dim M\neq 8\) we 
can take the limit \(R\to\infty\) and using the finiteness assumption \eqref{thm-finiteness-es4-harmonic-ass}
the calculation from above yields that 
\begin{align*}
\bar\Delta\tau(\phi)=0,\qquad R^N(d\phi(e_i),d\phi(e_j))\tau(\phi)=0. 
\end{align*}
The claim now follows by the same arguments given at the end of 
the proof of Theorem \ref{liouville-4-harmonic}.

\par\medskip
\textbf{Acknowledgements:}
The author would like to thank the reviewers for their many helpful comments
which helped to improve the presentation of the article's content.

The author gratefully acknowledges the support of the Austrian Science Fund (FWF) 
through the project P30749-N35 ``Geometric variational problems from string theory''.

\bibliographystyle{plain}
\bibliography{mybib}

\begin{thebibliography}{10}

\bibitem{MR655417}
P.~Baird and J.~Eells.
\newblock A conservation law for harmonic maps.
\newblock In {\em Geometry {S}ymposium, {U}trecht 1980 ({U}trecht, 1980)},
  volume 894 of {\em Lecture Notes in Math.}, pages 1--25. Springer, Berlin-New
  York, 1981.

\bibitem{MR2604617}
Paul Baird, Ali Fardoun, and Seddik Ouakkas.
\newblock Liouville-type theorems for biharmonic maps between {R}iemannian
  manifolds.
\newblock {\em Adv. Calc. Var.}, 3(1):49--68, 2010.

\bibitem{MR4106647}
V.~Branding, S.~Montaldo, C.~Oniciuc, and A.~Ratto.
\newblock Higher order energy functionals.
\newblock {\em Adv. Math.}, 370:107236, 60, 2020.

\bibitem{MR3834926}
Volker Branding.
\newblock A {L}iouville-type theorem for biharmonic maps between complete
  {R}iemannian manifolds with small energies.
\newblock {\em Arch. Math. (Basel)}, 111(3):329--336, 2018.

\bibitem{MR4007262}
Volker Branding.
\newblock The stress-energy tensor for polyharmonic maps.
\newblock {\em Nonlinear Anal.}, 190:111616, 17, 2020.

\bibitem{MR4184658}
Volker Branding.
\newblock A structure theorem for polyharmonic maps between {R}iemannian
  manifolds.
\newblock {\em J. Differential Equations}, 273:14--39, 2021.

\bibitem{MR4040175}
Volker Branding and Yong Luo.
\newblock A nonexistence theorem for proper biharmonic maps into general
  {R}iemannian manifolds.
\newblock {\em J. Geom. Phys.}, 148:103557, 9, 2020.

\bibitem{MR172310}
James Eells, Jr. and J.~H. Sampson.
\newblock \'{E}nergie et d\'{e}formations en g\'{e}om\'{e}trie
  diff\'{e}rentielle.
\newblock {\em Ann. Inst. Fourier (Grenoble)}, 14(fasc., fasc. 1):61--69, 1964.

\bibitem{MR2389639}
Fr\'{e}d\'{e}ric H\'{e}lein and John~C. Wood.
\newblock Harmonic maps.
\newblock In {\em Handbook of global analysis}, pages 417--491, 1213. Elsevier
  Sci. B. V., Amsterdam, 2008.

\bibitem{MR891928}
Guo~Ying Jiang.
\newblock The conservation law for {$2$}-harmonic maps between {R}iemannian
  manifolds.
\newblock {\em Acta Math. Sinica}, 30(2):220--225, 1987.

\bibitem{MR2395125}
E.~Loubeau, S.~Montaldo, and C.~Oniciuc.
\newblock The stress-energy tensor for biharmonic maps.
\newblock {\em Math. Z.}, 259(3):503--524, 2008.

\bibitem{MR3007953}
Shun Maeta.
\newblock The second variational formula of the {$k$}-energy and {$k$}-harmonic
  curves.
\newblock {\em Osaka J. Math.}, 49(4):1035--1063, 2012.

\bibitem{montaldo}
Stefano Montaldo and Alvaro Pampano.
\newblock Triharmonic curves in 3-dimensional homogeneous spaces.
\newblock {\em arXiv preprint arXiv:2008.10571}, 2020.

\bibitem{ou2019}
Ye-Lin Ou and Bang-Yen Chen.
\newblock {\em Biharmonic submanifolds and biharmonic maps in Riemannian
  geometry}.
\newblock World Scientific, 2019.

\bibitem{MR654088}
H.~C.~J. Sealey.
\newblock Some conditions ensuring the vanishing of harmonic differential forms
  with applications to harmonic maps and {Y}ang-{M}ills theory.
\newblock {\em Math. Proc. Cambridge Philos. Soc.}, 91(3):441--452, 1982.

\end{thebibliography}
\end{document}